\documentclass[11pt]{article}
\usepackage{amsthm}
\usepackage{epsf}
\usepackage{amsmath}
\usepackage{amsfonts}
\usepackage{amssymb}
\usepackage{amscd}
\usepackage[T1]{fontenc}
\usepackage{mathrsfs}
\usepackage{psboxit}

 \newtheorem{thm}{Theorem}[section]
 \newtheorem{cor}[thm]{Corollary}
 
 \newtheorem{prop}[thm]{Proposition}

 \newtheorem{rem}[thm]{Remark}

 \newcommand{\eps}{\varepsilon}
 \newcommand{\uep}{u^\eps}

 \newcommand{\R}{\mathbb R}
 \newcommand{\K}{\mathbb K}
 \newcommand{\N}{\mathbb N}
\newcommand{\lraup}{\relbar\joinrel\relbar\joinrel\relbar\joinrel\rightharpoonup}

 \newcommand{\tendef}[3]{{#1}\,\underset{#2}\lraup\,{#3}}
 \newcommand{\tende}[3]{\mbox{\raisebox{-3pt}{$\begin{CD}{#1}@>>{#2}>{#3}\end{CD}$}}}
 \newcommand{\fim}{\hfill{$\square$}}
 \newcommand{\cchi}{{\mbox{\raisebox{3pt}{$\chi$}}}}
 \newcommand{\bu}{\boldsymbol{u}}
\newcommand{\bU}{\boldsymbol{U}}
\newcommand{\bv}{\boldsymbol{v}}
\newcommand{\bff}{\boldsymbol{f}}
\newcommand{\bfF}{\boldsymbol{F}}
\newcommand{\bn}{\boldsymbol{n}}
\newcommand{\bw}{\boldsymbol{w}}
\newcommand{\bW}{\boldsymbol{W}}
\newcommand{\bk}{\boldsymbol{K}}
\newcommand{\bg}{\boldsymbol{g}}
\newcommand{\soma}{\displaystyle{\sum_{\mbox{\tiny $\begin{array}{c}1\!\!\le i_1\!\!\!<\!\!\ldots\!\!<i_k\!\!\!\le N\vspace{1mm}\\
i\in\{i_1,\ldots,i_k\}\end{array}$}}}}
\newcommand{\uniao}{\bigcup_{\mbox{\tiny $\begin{array}{c}1\!\!\le i_1\!\!\!<\!\!\ldots\!\!<i_k\!\!\!\le N\vspace{1mm}\\
i\in\{i_1,\ldots,i_k\}\end{array}$}}}

\newdimen\authorswidth
\newdimen\maxauthorswidth
\newcommand{\sameauthor}

\newcommand{\nextref}[4]{\bibitem{#1} #2, {\sl #3}, #4.
\global\setbox5=\hbox{#2}
\ifnum\wd5>\maxauthorswidth\authorswidth=\maxauthorswidth \else
\authorswidth=\wd5\fi
\renewcommand{\sameauthor}{\rule{\authorswidth}{0.5truept}}}

\title{\huge\bf On a constrained reaction-diffusion system related to multiphase problems}

\author{Jos\'{e} Francisco Rodrigues$\quad\qquad$ Lisa
Santos}

\date{}

\begin{document}

\maketitle


\vspace{5mm}

\begin{abstract}

We solve and characterize the Lagrange multipliers of a reaction-diffusion system in the Gibbs simplex of $\R^{N+1}$ by considering strong solutions of a system of parabolic variational inequalities in $\R^N$. Exploring properties of the two obstacles evolution problem, we obtain and approximate a $N$-system involving the characteristic functions of the saturated and/or degenerated phases in the nonlinear reaction terms. We also show continuous dependence results and we establish sufficient conditions of non-degeneracy for the stability of those phase subregions.

\end{abstract}

\baselineskip 16pt

\vspace{5mm}


\section{\bf Introduction} $\ $

This paper is motivated by the vector-valued reaction-diffusion
equation
\begin{equation}\label{system}
\partial_t\bU-\Delta \bU=\bfF(x,t,\bU),\qquad\mbox{ in }Q,
\end{equation}
for $\bU=\bU(x,t)$, defined from $Q=\,\Omega\times(0,T)$ into
$\R^{N+1}$, with homogeneous Neumann condition on
$\partial\Omega\times(0,T)$, where $\Omega$ is a bounded domain of
$\R^n$  and $T>0$ is arbitrary. We are interested in the case when
every component $u_i=u_i(x,t)$ is nonnegative and the system is
subject to the multiphase non-voids condition with
$\boldsymbol{J}=(1,\ldots,1)\in\R^{N+1}$:
\begin{equation}\label{2}
\bU\cdot \boldsymbol{J}=\sum_{j=1}^{N+1}u_j=1\qquad\mbox{ in }Q.
\end{equation}

From the equation (\ref{system}) it is clear that the constraint
(\ref{2}) implies $\bfF(x,t,\bU)\cdot \boldsymbol{J}=0$ in $Q$ and
so the reaction vector $\bfF$ should satisfy the necessary and very
restrictive condition
\begin{equation}\label{case}
F_{N+1}(x,t,V)=-\sum_{j=1}^N F_j(x,t,V)\ \ \mbox{ in }Q,\ \
\forall\,V=(v_1,\ldots,v_n,1-\sum_{j=1}^Nv_j),\ 0\le v_i\le 1.
\end{equation}

For instance, in replicator dynamics describing the evolution of
certain frequencies in a population, one possible definition of the
reaction term with this compatibility condition consists in choosing
\begin{equation}
F_{i}(x,t,V)= v_i
[\phi_i(x,t,V)-\sum_{j=i}^{N+1}v_j\phi_j(x,t,V)]\qquad \mbox{ in
}Q,\qquad i=1,...,N+1,
\end{equation}
where $v_i$ represents the $i$-frequency of the population and
$\phi_i$ the respective fitness (see, for instance, \cite{hs1} and \cite{hs2}), the constraint (\ref{2}) is essential to describe mixed strategies in evolutionary game theory in spatially homogeneous population dynamics (see \cite{m} and its references) or to model the non-voids condition in biological tissue growing \cite{k1,k2}. In phase fields models, the condition (\ref{2}) arises naturally in simulation of multiphase flows (\cite{kl}) and multiphase systems with diffuse phase boundaries, as in solidification of alloys or in grain boundary motion (see \cite{gsn} or \cite{bgn}).

Of course, in the case (\ref{case}), in particular, if $\bfF=\boldsymbol{0}$,
the problem becomes a simple one if the initial data $\bU(0)=\bU_0$
also satisfies the constraint (\ref{2}). However the situation is
entirely different in the general case of non trivial reactions,
specially in multiphase problems where at least one phase ``$i$'' in
a subregion of $Q$ is absent (i.e. $u_i=0)$, or fulfils another
subregion (when $u_i=1$).

Instead of solving the system (\ref{system}) in the Gibbs
(N+1)-simplex
$$\Psi = \{(v_1,\ldots,v_{N+1})\in R^{N+1}:\,\sum_{j=1}^{N+1}v_j= 1\mbox{ and }v_i\ge 0,\ i=1,\ldots, N+1\},$$
we shall replace this problem by the study of a unilateral problem
for the vector field of the first $N$ components
$\bu=(u_1,\ldots,u_N)$ of $\bU$, with the $N+1$ convex constraints
\begin{equation}\label{constraints}
\sum_{i=j}^Nu_j\le 1\qquad\mbox{ and }\qquad u_i\ge 0\qquad\mbox{ in
}Q,\qquad i=1,\ldots,N.
\end{equation}

This corresponds to solve the system of parabolic variational
inequalities, at each time $t\in(0,T)$,
\begin{multline}\label{vi}
\bu(t)\in\K:\qquad\int_\Omega\partial_t\bu(t)\cdot\big(\bv-\bu(t)\big)+\int_\Omega\nabla\bu(t)\cdot\nabla(\bv-\bu(t))\\
\ge\int_\Omega\bff(\bu(t))\cdot\big(\bv-\bu(t)\big),\qquad
\forall\bv\in\K,
\end{multline}
under the initial condition
\begin{equation}\label{ic}
\bu(0)=\bu_0=(u_{01},\ldots,u_{0N})\in\K.
\end{equation}

Here $\K$ denotes the convex subset of the Sobolev space
$H^1(\Omega)^N$ defined by
\begin{equation}\label{convex}
\K=\{\bv\in H^1(\Omega)^N:\, \sum_{j=1}^Nv_j\le 1,\ v_i\ge 0,\
i=1,\ldots,N,\ \mbox{in }\Omega\},
\end{equation}
where $\bv=(v_1,\ldots,v_N)$.

The reaction term $\bf$ may have a general form
$f_i(\bu)=f_i(x,t,\bU(x,t))$, $i=1,\ldots,N$, with $(x,t)\in\,Q$ and
$\bU=\big(u_1,\ldots,u_N,1-\displaystyle{\sum_{j=1}^Nu_j}\big)$. We
denote $\partial_t=\displaystyle{\frac{\partial\ }{\partial t}}$ and
$\nabla=\Big(\displaystyle{\frac{\partial\ }{\partial x_1}},\ldots,
\displaystyle{\frac{\partial\ }{\partial x_n}}\Big)$.

The main part of this work is the analysis of the new unilateral
problem (\ref{vi})-(\ref{ic}) under general assumptions on $\bff$:
only continuity on $\bu$ and integrability in $(x,t)\in\, Q$. In
particular, we prove that its solution $\bu=\bu(x,t)$, which each
component $u_i$ satisfies a double obstacle problem
\begin{equation}
0\le u_i\le 1-\sum_{j\neq i}u_j\qquad\mbox{ in }Q,\qquad
i=1,\ldots,N,
\end{equation}
where $\displaystyle{\sum_{j\neq i}u_j}$ denotes the sum of all
$N-1$ components but $u_i$ is, in fact, also the solution of a
reaction-diffusion system in the form
\begin{multline}\label{sistema}
\partial_t u_i-\Delta
u_i=f_i(\bu)+f_i^-(\bu)\cchi_{\{u_i=0\}}\\-\soma\frac1k\big(f_{i_1}(\bu)+\cdots+f_{i_k}(\bu)\big)^+\cchi_{i_i\ldots
i_k},\qquad\mbox{ in }Q.
\end{multline}

Here $\displaystyle{\soma}$ denotes the summation over all the
subsets $\{i_1,\ldots,i_k\}$  of $\{1,\ldots,N\}$ to which $i$
belongs, in particular, $k$ varies from $1$ to $N$. We also denote
 $g^+=g\vee 0$ and $g^-=-(g\wedge 0)$ the
positive and negative parts of a scalar function $g=g^+-g^-$,
$\cchi_A$ the characteristic function of the set $A$, (i.e.,
$\cchi_A=1$ in $A$ and $\cchi_A=0$ in $Q\setminus A$) and
$\cchi_{i_1\ldots i_k}$  the characteristic function of the set
$$I_{i_1\ldots i_k}=\big\{(x,t)\in\,Q:\big(u_{i_1}+\cdots+u_{i_k}\big)(x,t)=1,\ u_{i_j}(x,t)>0,\,
j=1,\ldots,k\},\quad k\in\{1,\ldots,N\}.$$

In particular $\displaystyle{\{u_i=1\}=\bigcap_{j\neq i}\{u_j=0\}}$,
i.e., one component is fully saturated if and only if the others are
absent. Hence from (\ref{sistema}) we see that, in general, the
respective reaction terms are coupled not only through the
semilinear term $\bff(\bu)$ but also through the characteristic
functions of the saturation sets of $I_{i_1\ldots i_k}$.

In this way, by setting for $i=1,\ldots,N,$
$$F_i(\bU)=f_i(\bu)+f_i^-(\bu)\cchi_{\{u_i=0\}}-\soma\frac1k\big(f_{i_1}(\bu)+\cdots+f_{i_k}(\bu)\big)^+\cchi_{i_i\ldots
i_k},$$
 with $\displaystyle{\bU=(\bu,1-\sum_{j=1}^N
u_j)}$, we can solve the system (\ref{system}) under the constraint
(\ref{2}) and identify the respective Lagrange multipliers
$h_i\equiv F_i(\bU)- f_i(\bU)$ in a precise form.

To illustrate the meaning of the system (\ref{sistema}), that contains $2^N-1+N$ characteristic functions, in general, we may consider the cases $N=1,2$ or $3$. Denoting, for simplicity, $f_i=f_i(\bu)$, $\cchi_i=\cchi_{\{u_i=1\}}$, we may write the Lagrange multipliers as
\begin{equation*}
\begin{array}{l}
h_1=f_1^-\cchi_{\{u_1=0\}}-f_1^+\cchi_1-\frac12(f_1+f_2)^+\cchi_{12}-\frac12(f_1+f_3)^+\cchi_{13}-\frac13(f_1+f_2+f_3)^+\cchi_{123}
\vspace{3mm}\\
h_2=f_2^-\cchi_{\{u_2=0\}}-f_2^+\cchi_2-\frac12(f_1+f_2)^+\cchi_{12}-\frac12(f_2+f_3)^+\cchi_{23}-\frac13(f_1+f_2+f_3)^+\cchi_{123}
\vspace{3mm}\\
h_3=f_3^-\cchi_{\{u_3=0\}}-f_3^+\cchi_3-\frac12(f_1+f_3)^+\cchi_{13}-\frac12(f_2+f_3)^+\cchi_{23}-\frac13(f_1+f_2+f_3)^+\cchi_{123}
\end{array}
\end{equation*}

Ignoring the third equation and all the terms involving the third component, we may obtain the case $N=2$. The first two terms of the right hand side of the first equation correspond, in the case $N=1$, to the scalar two obstacles problem that has been proposed for phase separations in \cite{be1,be2}.

The mathematical treatment of this unilateral system is done in the following three sections. In section 2, we consider the semilinear approximation of the unique solution of (\ref{vi})-(\ref{ic}) in the case of the reaction $\bff$ is in $L^2(Q)^N$ and independent of the solution. Although there exists a large literature on parabolic variational inequalities (see, for instance, \cite{jll}, \cite{b72}, \cite{ken75}, \cite{ct78} or \cite{f82}), the direct approach of the  bounded penalization used for the two obstacles problem in \cite{r2} (see also \cite{p-00}), extended here for the system (\ref{sistema}), allows the use of monotone methods. This yields a direct way of obtaining Lewy-Stamppachia inequalities (\ref{pui}), obtained first by \cite{ct78} for parabolic problems, implying the $W^{2,1}_p$ and H\"{o}lder regularity for the solution to (\ref{vi}). Similar results for the $N$-membranes stationary problem have been obtained in \cite{ars1,ars2}. We note in our case the simplification due to homogeneous Neumann condition.

In section 3, we extend the existence result to general nonlinear reaction $\bff=\bff(\bu)$ taking values in $L^1(Q)^N$. Here we explore the fact that the convex set (\ref{convex})  lies in the unit disc and we extend the direct technique of \cite{pu}. We show also a continuous dependence result and, in the case of $\lambda I-\bff$ being monotone non-decreasing, in particular if $\bff$ is Lipschitz continuous in $\bu$, also the uniqueness of solution and their strong approximation by the penalized solutions.

Finally, in the last section, we characterize the solution of the variational inequality (\ref{vi}) as solutions of the reaction-diffusion system (\ref{sistema}), by extending  some remarks of \cite{r3} to the two obstacles parabolic problem. We also show that
$$\{u_i=0\}\subset\{f_i(\bu)\le 0\}\qquad\mbox{ and }\qquad I_{i_1\ldots i_k}\subset\big\{\sum_{j=1}^kf_{i_j}(\bu)\ge 0\big\}$$
a.e. in $Q$, for $1\le i_1<\cdots<i_k\le N$, $\forall\,k=1,\ldots,N$ and we can modify the system (\ref{sistema}) (see (\ref{eq70})) and show that the a.e. pointwise nondegeneracy assumptions
$$\sum_{j=1}^kf_{i_j}(\bu)\neq 0,\qquad 1\le i_1<\cdots<i_k\le N,\quad
k=1,\ldots,N,$$
are sufficient conditions for the local stability of the characteristic functions $\cchi_{\{u_i=0\}}$ and $\cchi_{i_1\ldots i_k}$ with respect  to the perturbation of the nonlinear reaction terms $\bff$.

\vspace{5mm}

\section{Approximation of strong solutions by semilinear problems}$ $

In this section we consider the case where $\bff=(f_1,\ldots,f_N)$
depends only on $(x,t)$ and is given in $L^2(Q)^N$.

To prove existence of solution of the variational inequality
(\ref{vi})-(\ref{ic}), we consider a family of approximating
semilinear systems of equations. We define, for each $\eps>0$,
$\theta_\eps:\R\longrightarrow\R$ by
\begin{equation}\label{thetaeps}
\theta_\eps(s)=\left\{\begin{array}{ll} 0&\mbox{ if }s\ge
0\vspace{2mm}\\
s/\eps&\mbox{ if }-\eps< s< 0\vspace{2mm}\\
-1 &\mbox{ if }s\le-\eps,
\end{array}
\right.
\end{equation}
and we denote
$$P\bu=\partial_t\bu-\Delta\bu=(Pu_1,\ldots,Pu_N),$$ where
$\partial_t\bu=(\partial_tu_1,\ldots,\partial_tu_N)$ and
$\Delta\bu=(\Delta u_1,\ldots,\Delta u_N)$. We also denote
$Pu_i=\partial_tu_i-\Delta u_i$, $i=1,\ldots,N$. The approximating
problems are given by the following weakly coupled parabolic system
with Neumann condition
\begin{alignat}{1}\label{aprox}
&P\uep_i+f_i^-\theta_\eps(\uep_i)-\soma\frac1k\big(f_{i_1}+\cdots+f_{i_k}\big)^+\theta_\eps(1-\uep_{i_1\ldots i_k})=f_i\qquad\mbox{ in }Q,\\
\label{lateral}&\frac{\partial\uep_i}{\partial\bn}=0\qquad\mbox{ on
}\partial\Omega\times(0,T), \\
\label{inicial} &\uep_i(0)=u_{0i}\qquad\mbox{ in
}\Omega,\qquad\qquad\qquad\qquad(i=1,\ldots,N)
\end{alignat}
where $\displaystyle{\frac{\partial\ }{\partial\bn}}$ is the outward normal
derivative on $\partial\Omega\times(0,T)$, the meaning of $\soma$
was explained in the introduction and
\begin{equation}\label{v1i}
\forall\,\bv=(v_1,\ldots,v_N)\ \
\forall\,\{i_1,\ldots,i_k\}\subseteq\{1,\ldots,N\}\qquad
v_{i_1\ldots i_k}=v_{i_1}+\cdots+v_{i_k}.
\end{equation}

Defining the penalization operator $\Theta_\eps$ by
\begin{alignat}{1}\label{thetaeps2}
\Theta_\eps\bu\cdot\bv=&\sum_{i=1}^N\Big[f_i^-\theta_\eps(u_i)-
\soma\frac1k\big(f_{i_1}+\cdots+f_{i_k}\big)^+\theta_\eps(1-u_{i_1\ldots
i_k})\Big]v_i\\
=&\sum_{i=1}^Nf_i^-\theta_\eps(u_i)v_i-\sum_{\mbox{\tiny $1\!\!\le
i_1\!\!\!<\!\!\ldots\!\!<i_k\!\!\!\le
N$}}\frac1k\big(f_{i_1}+\cdots+f_{i_k}\big)^+\theta_\eps(1-u_{i_1\ldots
i_k})v_{i_1\ldots i_k},
\end{alignat}
 we formulate
(\ref{aprox})-(\ref{lateral}) in variational form for a.e.
$t\in(0,T),$
\begin{equation}\label{iveps}
\int_\Omega\partial_t\bu^\eps(t)\cdot\bv+\int_\Omega\nabla\bu^\eps(t)\cdot\nabla\bv+\int_\Omega\Theta_\eps(\bu^\eps(t))\cdot\bv
=\int_\Omega\bff(t)\cdot \bv,\qquad\forall\bv\in\, H^1(\Omega)^N,
\end{equation}
associated with the initial condition (\ref{inicial}).

\begin{prop} \label{est} Assuming that
\begin{equation}\label{fu0}
\bff=(f_1,\ldots,f_N)\in L^2(Q) ^N\qquad\mbox{ and }\qquad
\bu_0\in\K,
\end{equation}
the problem {\em(\ref{iveps})-(\ref{inicial})} has a unique solution
$\bu^\eps\in
  H^1(0,T;L^2(\Omega)^N)\cap\, L^\infty(0,T;H^1(\Omega)^N))$.
\end{prop}
\proof We begin by proving the monotonicity of the penalization
operator $\Theta_\eps$.

In fact, recalling that $\theta_\eps$ is monotone nondecreasing and
 the definition (\ref{v1i}) we have
\begin{alignat*}{1}
\big(\Theta_\eps\bu-\Theta_\eps\bv\big)&\cdot\big(\bu-\bv\big)
\\
&\hspace{-15mm}=\sum_{i=1}^Nf_i^-\big(\theta_\eps(u_i)-\theta_\eps(v_i)\big)(u_i-v_i)\\
&\hspace{-15mm}-\sum_{\mbox{\tiny $1\!\!\le i_1\!\!\!<\!\!\ldots\!\!<i_k\!\!\!\le
N$}}\frac1k\big(f_{i_1}+\cdots f_{i_k}\big)^+\big(\theta_\eps(1-u_{i_1\ldots i_k})-\theta_\eps(1-v_{i_1\ldots i_k})\big)(u_{i_1\ldots i_k}-v_{i_1\ldots i_k}),\\
&\hspace{-15mm}\ge 0,
\end{alignat*}
since $f_j^-$ and $\big(f_{i_1}+\cdots f_{i_k}\big)^+$ are
nonnegative functions.

The existence and uniqueness of solution  $\bu^\eps\in L^2(0,T;
H^1(\Omega) ^N)$ is immediate by applying the theory of monotone
operators (\cite{jll}, \cite{zheng})).

Setting $\bv=(\uep_1,\ldots,\uep_N)$ in the approximating problem
(\ref{iveps}) and integrating in time, letting
$$g^\eps_i=P\uep_i=f_i-f_i^-\theta_\eps(\uep_i)+\soma\frac1k\big(f_{i_1}+\cdots+f_{i_k}\big)^+\theta_\eps(1-\uep_{i_1\ldots
i_k}),$$ which is bounded in $L^2(Q)$ independently of $\eps$, we
obtain that, for every $0<t<T$, with $Q_t=\Omega\times(0,t)$,
$$\frac12\int_\Omega|\bu^\eps(t)|^2+\int_{Q_t}|\nabla\bu^\eps|^2\le\frac12\int_\Omega
|\bu_0|^2+\frac12
\int_{Q_t}|\boldsymbol{g}^\eps|^2+\frac12\int_{Q_t}|\bu^\eps|^2.$$

The Grownwall inequality yields the uniform boundedeness (in $\eps$)
of $\bu^\eps$, first in\linebreak
$L^\infty(0,T;L^2(\Omega)^N))$ and afterwards also in
$L^2(0,T; H^1(\Omega) ^N)$.

Letting, formally, $\bv=\partial_t\bu^\eps$ in (\ref{iveps}) (in
fact in the respective Faedo-Galerkin approximation) and integrating
in time, we get
$$\int_{Q_t}\big|\partial_t\bu^\eps\big|^2+\int_\Omega|\nabla\bu^\eps(t)|^2\le\int_{Q_t}|\boldsymbol{g}^\eps|^2+\int_\Omega|\nabla
\bu_0|^2$$ and so $\partial_t\bu_\eps$ is also bounded in $L^2(Q)
^N$ and $\nabla\bu^\eps$ in $L^\infty(0,T;L^2(\Omega)^N)$. Therefore
\begin{equation}
\label{esteps2} \{\bu^\eps\}_{\eps>0}\ \mbox{ is bounded in }\
H^1(0,T;L^2(\Omega)^N))\cap\,L^\infty(0,T;H^1(\Omega)^N).\end{equation}
 \fim\vspace{2mm}

\begin{prop}  \label{desigeps} Assuming {\em (\ref{fu0})}, the solution $\bu^\eps$ of the problem
{\em(\ref{iveps})-(\ref{inicial})} satisfies
\begin{equation}\label{convexo}
\uep_i\ge-\eps,\qquad i=1,\ldots,N,\qquad\qquad\sum_{i=1}^N\uep_i\le
1+\eps.
\end{equation}
\end{prop}
\proof In fact, we are going to prove the following more general set
of inequalities
\begin{equation*}
\uep_i\ge-\eps,\quad i=1,\ldots,N,\qquad\mbox{ and
}\qquad\uep_{i_1\ldots i_r}\le 1+\eps,\quad \forall\,1\le
i_1<\ldots<i_r\le N
\end{equation*}
and the proof of the right hand side inequalities will be done by
induction on $r$.

Let us prove the case $r=1$, i.e., $\uep_i\le 1+\eps$, for all
$i\in\{1,\ldots,N\}$. Multiplying the $i$-th equation of the
approximating system (\ref{aprox}) by $(\uep_i-(1+\eps))^+$ and
integrating over $Q_t=\Omega\times(0,t)$, we have
\begin{multline*}
\int_{Q_t}\partial_t\uep_i(\uep_i-(1+\eps))^++\int_{Q_t}\nabla\uep_i\cdot\nabla(\uep_i-(1+\eps))^+\\
=\int_{Q_t}\Big[f_i-f_i^-\theta_\eps(\uep_i)+
\soma\frac1k\big(f_{i_1}+\cdots+f_{i_k}\big)^+\theta_\eps(1-\uep_{i_1\ldots
i_k})\Big](\uep_i-(1+\eps))^+
\end{multline*}

Recalling that $-1\le\theta_\eps\le 0$ and that, in the set
$\{\uep_i>1+\eps\}$, we have $\theta_\eps(\uep_i)=0$ and
$\theta_\eps(1-\uep_i)=-1$ we get
\begin{equation}\label{111}
\frac12\int_\Omega|(\uep_i-(1+\eps))^+(t)|^2+\int_{Q}|\nabla(\uep_i-(1+\eps))^+|^2\le\int_Q
(f_i-f^+_i)(\uep_i-(1+\eps))^+\le 0,
\end{equation}
so $(\uep_i-(1+\eps))^+\equiv 0$, i.e. $\uep_i\le 1+\eps$.

Assuming we have proved that $\uep_{i_1\ldots i_r}\le 1+\eps$, we
are going to show that $\uep_{i_1\ldots i_r i_{r+1}}\le 1+\eps$.

We multiply the   equations $i_j$, $j=1,\ldots,r+1$, by
$(\uep_{i_1\ldots i_ri_{r+1}}-(1+\eps))^+$, sum from $1$ to $r+1$
and integrate over $Q_t$. We obtain
\begin{multline*}
\int_{Q_t} P\uep_{i_1\ldots i_ri_{r+1}}(\uep_{i_1\ldots
i_ri_{r+1}}-(1+\eps))^+=\int_{Q_t}\Big[\sum_{j=1}^{r+1}f_{i_j}-\sum_{j=1}^{r+1}f_{i_j}^-
\theta_\eps(\uep_{i_j})\\
+\sum_{j=1}^{r+1}
\sum_{\mbox{\tiny $\begin{array}{c}1\!\!\le i_1\!\!\!<\!\!\ldots\!\!<i_k\!\!\!\le N\vspace{1mm}\\
i_j\in\{i_1,\ldots,i_k\}\end{array}$}}\frac1k\big(f_{i_1}+\cdots+f_{i_k}\big)^+\theta_\eps(1-\uep_{i_1\ldots
i_k})\Big](\uep_{i_1\ldots i_ri_{r+1}}-(1+\eps))^+.
\end{multline*}

Observe that, in the set $\{\uep_{i_1\ldots i_ri_{r+1}}>1+\eps\}$ we
have $\uep_{i_j}\ge 0$, for $j=1,\ldots,r+1$, since, by induction,
$\uep_{l_1\ldots l_r}=\uep_{i_1}+\cdots+\uep_{i_{r+1}}-\uep_{i_j}\le 1+\eps$. So, in that set
$\theta_\eps(\uep_{i_j})=0$ and, on the other hand,
$\theta_\eps(1-\uep_{i_1\ldots i_ri_{r+1}})=-1$. The induction conclusion follows from
\begin{multline*}
\int_{\Omega}
|(\uep_{i_1\ldots i_ri_{r+1}}-(1+\eps))^+(t)|^2+\int_{Q_t}|\nabla(\uep_{i_1\ldots i_ri_{r+1}}-(1+\eps))^+|^2\\
\le\int_{Q_t}\Big[\sum_{j=1}^{r+1}f_{i_j}-(r+1)
\frac1{r+1}\big(f_{i_1}+\cdots+f_{i_{r+1}}\big)^+\Big](\uep_{i_1\ldots
i_ri_{r+1}}-(1+\eps))^+\le 0.
\end{multline*}

To prove that $\uep_i\ge -\eps$, we multiply  the $i$-th equation of
(\ref{aprox}) by $(-\uep_i-\eps)^+$, obtaining
\begin{multline*}
\frac12\int_\Omega|(-\uep_i-\eps)^+(t)|^2+\int_{Q}|\nabla(-\uep_i-\eps)^+|^2=\int_Q\Big[-f_i+f_i^-\theta_\eps(\uep_i)\\
-\soma\frac1k\big(f_{i_1}+\cdots+f_{i_k}\big)^+\theta_\eps(1-\uep_{i_1\ldots
i_k})\Big](-\uep_i-\eps)^+.
\end{multline*}

Let $J_{k,i}=\{i_1,\ldots,i_k\}\setminus\{i\}$ and denote the
elements of $J_{k,i}$ by $j_1\ldots j_{k-1}$. Since, in the set
$\{(-\uep_i-\eps)^+>0\}=\{\uep_i<-\eps\}$, we have
$1-\uep_{i_1\ldots i_k}=1-\uep_{j_1\ldots j_{k-1}}-\uep_i>0$ (recall
that $\uep_{j_1\ldots j_{k-1}}\le 1+\eps$). So,
\begin{equation*}
\frac12\int_\Omega|(-\uep_i-\eps)^+(t)|^2+\int_{Q}|\nabla(-\uep_i-\eps)^+|^2\le\int_Q\Big[-f_i-f_i^-\Big](-\uep_i-\eps)^+
=\int_Q-f_i^+(-\uep_i-\eps)^+\le 0,
\end{equation*}
that implies  $(-\uep_i-\eps)^+=0$, or $\uep_i\ge-\eps$.
 \fim\vspace{2mm}

 \begin{thm}  \label{existence} Assuming {\em (\ref{fu0})}, the variational inequality {\em
(\ref{vi})-(\ref{ic})} has a unique solution $\bu$ such that
\begin{equation}\label{23}
\bu\in\,H^1(0,T;L^2(\Omega)^N)\cap\, L^\infty(0,T;H^1(\Omega)^N)
\end{equation}
and
\begin{equation}\label{24}
P\bu\in\,L^2(Q)^N.
\end{equation}
\end{thm}
\proof Let $\bu^\eps$ be the solution of the problem (\ref{iveps}).
Using the uniform estimates (in $\eps$) obtained in (\ref{esteps2}),
we know  there exists $\bu$ such that
\begin{equation*}\label{conv}
\begin{array}{l}\tende{\bu^\eps}{\eps}{\bu} \qquad\mbox{ in
} L^2(Q) ^N\mbox{ strong},\vspace{3mm}\\
\tendef{\bu^\eps}{\eps}{\bu}\qquad\mbox{ in
}L^\infty(0,T; H^1(\Omega) ^N)\mbox{ weak-}*,\vspace{3mm}\\
\tendef{\partial_t\bu^\eps}{\eps}{\partial_t\bu}\qquad\mbox{ and
}\qquad P\bu^\eps\lraup P\bu\qquad\mbox{ in } L^2(Q) ^N\mbox{ weak}.
\end{array}
\end{equation*}

We have $\bu(t)\in\K$, for a.e. $t\in[0,T]$, because $\bu^\eps$
satisfies the inequalities (\ref{convexo}).

Given $\bv\in\, L^2(0,T;\K)$, set $\bv(t)-\bu^\eps(t)$ in
(\ref{iveps}) and integrate in time. Then
$$\int_Q\partial_t\bu^\eps\cdot
(\bv-\bu)+\int_Q\nabla\bu^\eps\cdot\nabla (\bv-\bu)\ge \int_Q
\bff^\eps\cdot(\bv-\bu^\eps),$$ since
$\displaystyle{\int_Q\big(\Theta_\eps(\bu^\eps)-\Theta_\eps(\bv)\big)\cdot(\bv-\bu^\eps)\le
0}$ and $\Theta_\eps(\bv(t))=0$ if $\bv(t)\in\K$. Passing to the
limit when $\eps\rightarrow 0$ and noting that
$$\liminf_{\eps\rightarrow
0}\int_Q\big(\partial_t\bu^\eps\cdot\bu^\eps+\nabla\bu^\eps\cdot\nabla\bu^\eps\big)\ge\int_Q\big(\partial_t\bu\cdot\bu+\nabla\bu\cdot\nabla\bu\big),$$
we find that $\bu$ satisfies (\ref{ic}) and
\begin{equation}\label{vit}
\int_Q\partial_t \bu\cdot(\bv-\bu)+\int_Q\nabla \bu\cdot\nabla
(\bv-\bu)\ge \int_Q \bff\cdot(\bv-\bu),\qquad\forall\bv\in\,
L^2(0,T;\K),\end{equation}
 which is easily seen to be equivalent to
(\ref{vi}). The uniqueness is immediate.
 \fim\vspace{2mm}

 We remark that no regularity of the boundary $\partial\Omega$ has
 been required in (\ref{iveps}) and, in fact, the Neumann boundary
 condition (\ref{lateral}) is only formal.  In the proof of Theorem \ref{existence}
 we have used the compactness of the sequence $\{\bu^\eps\}_\eps$
 in $L^2(Q)^N$. This holds, for instance, for domains with Lipschitz boundaries, but also, since the sequence $\{\bu^\eps\}_\eps$
 is uniformly bounded in $L^\infty(Q)^N$,  for a larger class of bounded open subsets of $\R^{N+1}$.
 However, the
 approximation by semilinear parabolic equations yields immediately
 an additional regularity of these strong solutions.

 Indeed, from the definitions of $\theta_\eps$ and $\Theta_\eps$, from (\ref{iveps}) with arbitrary $\varphi\in\mathcal{D}(Q)$,
 $\varphi\ge 0$, we find
 \begin{equation}\label{pui}
 f_i-\soma\frac1k\big(f_{i_1}+\cdots+f_{i_k}\big)^+\le P\uep_i=f_i-\Theta_\eps(\bu^\eps)\le
 f_i+f_i^-=f_i^+\qquad\mbox{ a.e. in }Q.
 \end{equation}

 By the conclusion of Theorem \ref{existence} we also obtain,
 for each $i=1,\ldots,N,$
 \begin{equation}\label{puiu}
 f_i-\soma\frac1k(f_{i_1}+\cdots+f_{i_k})^+\le Pu_i\le
f^+_i\qquad\mbox{ a.e. in }Q
 \end{equation}
 and we can apply directly the second order linear parabolic theory
 (see \cite{lsu}) in the Sobolev spaces
 $$W^{2,1}_p(Q)=W^{1,p}(0,T;L^p(\Omega))\,\cap\,
 L^p(0,T;W^{2,p}(\Omega)),\qquad 1<p<\infty.$$

 These spaces satisfy the Sobolev imbeddings, for $p>(n+2)/(2-k)$,
 with $k=0,1,$
 $$W^{2,1}_p(Q)\,\subset\, C^{k,0}_\alpha(\overline{Q}),\qquad
 0\le\alpha<2-k-(n+2)/p,$$
 where $C^{k,0}_\alpha(\overline{Q})$ denotes the spaces of
 H\"{o}lder continuous functions $v$ in $Q$, with exponent $\alpha$
 in the $x$-variables and $\alpha/2$ in the $t$-variable and, in the
 case $k=1$, with $\nabla v$ satisfying the same property (see
 \cite{lsu}, p. 80). Therefore, as a consequence of (\ref{puiu}), we conclude. \fim\vspace{2mm}

 \begin{thm} \label{t24}Assume that $\partial\Omega$ is smooth, say of class
 $C^2$ and
 \begin{equation}
 \bff\in\, L^p(Q)^N\qquad\mbox{ and }\qquad \bu_0\in\,\K\cap
 W^{2-2/p,p}(\Omega)^N,\qquad 1<p<\infty,
 \end{equation}
 with each component $u_{0i}$ satisfying the compatibility condition
 $\displaystyle{\frac{\partial u_{0i}}{\partial\boldsymbol{n}}=0}$ on $\partial\Omega$ if
 $p>3$.

 \vspace{1mm}

 Then the unique solution $\bu$ of the variational inequality
 {\em (\ref{vi})-(\ref{ic})} is such that
 \begin{equation}
 \bu\in\, W^{2,1}_p(Q)^N\,\cap\, L^\infty(0,T;\K),
 \end{equation}
 and, in particular, is H\"{o}lder continuous in $\overline{Q}$ if
 $p>(n+2)/2$ and has $\nabla\bu$ also H\"{o}lder continuous if
 $p>n+2$.
 \end{thm}\fim\vspace{2mm}

We observe that, when $p<2$, the inclusion $W^{2,1}_p(Q)\,\subset\,
L^2(0,T;H^1(\Omega))$ only takes place if $p\ge(2n+4)/(n+4)$ but, as
we shall see in the next section and since $\K$ is bounded,
(\ref{vi})-(\ref{ic}) is solvable for any $\bff\in\, L^1(Q)^N$.

\vspace{5mm}

\section{Existence and uniqueness of variational solutions}$ $

In this section, requiring the compactness of the inclusion of
$H^1(\Omega)$ into $L^2(\Omega)$ by assuming a Lipschitz boundary
$\partial\Omega$, we show how we can still solve the variational
inequality (\ref{vit}) for a more general initial condition
\begin{equation}\label{ktilde}
\bu_0\in\,\tilde{\K}=\{v\in\, L^2(\Omega)^N:\,\sum_{j=1}^Nv_j\le 1,\
v_i\ge 0,\ i=1,\ldots,N,\mbox{ in }\Omega\}\end{equation} and for
general nonlinear $\bff=\bff(\bu)$ defining a continuous operator
from $L^2(0,T;\tilde{\K})$ in $L^1(Q)^N$. We shall assume that
$\bff=\bff(x,t,\bv):Q\times\,[0,1]^N\rightarrow\R^N$ satisfies
\begin{equation}\label{fnova1}
\bff=\bff(x,t,\bv)\mbox{ is continuous in }\bv\mbox{ for a.e.
}(x,t)\in\,Q,
\end{equation}
\begin{equation}\label{fnova2}
\exists\,\varphi_1\in\,L^1(Q):\qquad|\bff(x,t,\bv)|\le\varphi_1(x)\qquad\forall\,\bv\in\,[0,1]^N,\
\mbox{ for a.e. }(x,t)\in\,Q.
\end{equation}

However, now the solution has less regularity, namely
\begin{equation}\label{regu}
\bu\in\,C([0,T];L^2(\Omega)^N\cap\,\tilde{\K})\,\cap\,
L^2(0,T;H^1(\Omega)^N)
\end{equation}
and its derivative may not be a function, since we only have
\begin{equation}\label{regut}
\partial_t\bu\in\,
L^1(Q)^N\,+\,L^2\big(0,T;\big(H^1(\Omega)^N\big)'\big).
\end{equation}

Hence the first term in the variational inequality (\ref{vit})
should be interpreted  in the duality sense between
$L^1(Q)^N\,+\,L^2\big(0,T;\big(H^1(\Omega)^N\big)'\big)$ and
$L^\infty(Q)^N\,\cap\,L^2(0,T;H^1(\Omega)^N)$, namely through the
formula
\begin{equation}\label{langle}
\langle\partial_t\bu,\bv\rangle_t=\int_{Q_t}P\bu\cdot\bv-\int_{Q_t}\nabla\bu\cdot\nabla\bv,\qquad\forall\,\bv\in\,
L^\infty(Q)^N\,\cap\,L^2(0,T;H^1(\Omega)^N),
\end{equation}
for arbitrary $t\in\,(0,T]$ since, as we shall see, (\ref{puiu})
yields $P\bu\in\,L^1(Q)^N$.

\begin{thm} \label{fu}
Under the assumptions {\em (\ref{ktilde})}, {\em (\ref{fnova1})}and
{\em (\ref{fnova2})}, the variational inequality {\em (\ref{vit})}
has a solution $\bu$ satisfying {\em (\ref{regu})}, {\em
(\ref{regut})}, {\em (\ref{puiu})} and $\bu(0)=\bu_0$ and we can
write
\begin{equation}\label{vi2}
\int_Q\big(P\bu-\bff(\bu)\big)\cdot(\bv-\bu)\ge
0,\qquad\forall\,\bv\in\, L^2(0,T;\tilde{\K}).
\end{equation}
\end{thm}
\proof We consider the closed convex subset of $L^2(Q)^N$
$$\bk=L^2(0,T;\tilde{\K})=\{\bv\in\,L^2(Q)^N:\, u_i\ge 0,\,
i=1,\ldots,N,\ \sum_{i=1}^Nu_i\le 1\mbox{ in }Q\}$$ and we define
$\Phi:\bk\rightarrow\bk$ as the nonlinear operator that associates
to each $\bw\in\,\bk$ the solution $\bu_{\bw}=\Phi(\bw)$ of the
variational inequality (\ref{vit}) with $\bff$ replaced by
$\boldsymbol{g}=\bff(x,t,\bw)$ and fixed initial data
$\bu_0\in\,\tilde{\K}$.

By showing that $\Phi$ is a continuous and compact operator, a fixed
point $\bu=\Phi(\bu)$, given by Schauder Theorem, will provide a
solution with the required properties.

Indeed, first we observe that if we consider any sequence $\bk\ni
\tende{\bw_\nu}{\nu}{\bw}\in\bk$ in $L^2(Q)^N$, by (\ref{fnova1})
and (\ref{fnova2}), the Lebesgue Theorem implies
$$\tende{\bg_\nu=\bff(\bw_\nu)}{\nu}{\bff(\bw)=\bg}\qquad\mbox{ in } L^1(Q)^N.$$

Next, for any $\bg\in\,L^1(Q)^N$ and any $\bu_0\in\tilde{\K}$ we
consider sequences  $\bg_\nu\in\,L^2(Q)^N$ and $\bu_{0\nu}\in\K$
such that
$$\tende{\bg_\nu}{\nu}{\bg}\ \mbox{ in } L^1(Q)^N\qquad\mbox{ and
}\qquad\tende{\bu_{0\nu}}{\nu}{\bu_0}\ \mbox{ in } L^2(\Omega)^N$$
and we denote by $\bu_\nu\equiv S(\bu_{0\nu},\bg_\nu)$ the unique
solution of (\ref{vit})-(\ref{ic}) given by Theorem \ref{existence},
for each $\bg_\nu$ and $\bu_{0\nu}$. We observe that each component
of $P\bu_\nu$ satisfies the inequality (\ref{puiu}) with $f_i$
replaced by $(\bg_\nu)_i$. From (\ref{vit}) for $\bu_\mu$ and
$\bu_\nu$, we easily find, for a.e. $t\in(0,T)$
$$\frac12\frac{d\ }{dt}\int_\Omega|\bu_\mu-\bu_\nu|^2+\int_\Omega|
\nabla(\bu_\mu-\bu_\nu)|^2\le\int_\Omega(\bg_\mu-\bg_\nu)\cdot(\bu_\mu-\bu_\nu)$$
and, integrating in time, we obtain
\begin{equation}\label{cauchy}
\sup_{0<t<T}\int_\Omega|\bu_\mu(t)-\bu_\nu(t)|^2+\int_Q|
\nabla(\bu_\mu-\bu_\nu)|^2\le\int_\Omega|\bu_{0\mu}-\bu_{0\nu}|^2+4\int_Q|\bg_\mu-\bg_\nu|.
\end{equation}

This estimate shows that $\{\bu_\nu\}_\nu$ is a Cauchy sequence in the
Banach space
\begin{equation}
\bW=C([0,T];L^2(\Omega)^N)\,\cap\, L^2(0,T;H^1(\Omega)^N)
\end{equation}
with respect to the norm
\begin{equation}\label{norm}
\|\hspace{-1pt}|\bv\|\hspace{-1pt}|=\Big(\sup_{0<t<T}\int_\Omega|\bv(t)|^2+\int_Q|\nabla\bv|^2\Big)^{1/2}
\end{equation}
and, hence, there exists a function $\bu_{\bg}\in\,\bW$
$$\tende{\bu_\nu}{\nu}{\bu_{\bg}}\qquad\mbox{ in }\bW.$$

In addition, $\bu_{\bg}\in L^2(0,T;\K)\,\cap\,C([0,T];\tilde{\K})$
and $P\bu_{\bg}\in\, L^1(Q)^N$, which implies, by (\ref{langle}),
that $\partial_t\bu_{\bg}$ satisfies (\ref{regut}). Hence, using
(\ref{langle}), we may pass to the limit in $\nu$ in
$$\langle P\bu_\nu-\bg_\nu,\bv-\bu_\nu\rangle=\int_Q(P\bu_\nu-\bg_\nu)\cdot(\bv-\bu_\nu)\ge
0$$ for an arbitrary $\bv\in L^2(0,T;\K)\subset L^\infty(Q)^N$, and
using the formula
$$2\langle\partial_t\bu_{\bg},\bu_{\bg}\rangle_t=\int_\Omega|\bu_{\bg}(t)|^2-\int_\Omega|\bu_0|^2,\qquad\forall\,t\in\,(0,T],$$
we conclude that $\bu_{\bg}=S(\bu_0,\bg)$ is the (unique) solution
of the variational inequality (\ref{vit}) (or equivalently
(\ref{vi2})) with data $\bg\in\,L^1(Q)^N$ and
$\bu_0\in\,\tilde{\K}$. In particular, from (\ref{cauchy}), we also
obtain that, for fixed $\bu_0\in\,\tilde{\K}$, the operator
$\Sigma:\bg\mapsto\bu_{\bg}=S(\bu_0,\bg)$ is H\"{o}lder continuous of
order $1/2$, from $L^1(Q)^N$ into $\bW$.

Since $\partial_t\bu_{\bg}$ satisfies the property (\ref{regut}), it
is in fact in $L^1(0,T;H^{-s}(\Omega)^N)$, for $s$ sufficiently
large and, by a well known compactness embedding (see \cite{si} or
Theorem 3.11 of \cite{zheng}), the compactness of
$H^1(\Omega)\subset L^2(\Omega)$ implies that, in fact, $\Sigma$
regarded as an operator from $L^1(Q)^N$ into $\bk\subset L^2(Q)^N$
is, therefore, completely continuous. Hence, $\Phi=\Sigma\circ\bff$
fulfils the requirements of the Schauder fixed point theorem and the
proof is complete.\fim\vspace{2mm}

\begin{rem} It is clear that if $\bu_0\in\K$ and, in {\em
(\ref{fnova2})}, $\varphi_1\in L^2(Q)$, we obtain in {\em Theorem
\ref{fu}} the existence of a strong solution satisfying {\em
(\ref{23})} and {\em (\ref{24})}. Of course, if we have the
regularity assumptions of {\em Theorem \ref{t24}}, i.e.,
$\varphi_1\in L^p(Q)$, implying by the inequalities {\em
(\ref{puiu})} that $P\bu\in L^p(Q)^N$, we also obtain solutions in
$W^{2,1}_p(Q)^N$, in particular H\"{o}lder continuous solutions if
$p>(n+2)/2$.
\end{rem}

In general (\ref{vi2}) may have more than one solution, but if we
assume, in addition, that for some $\lambda>0$, $\lambda\,I-\bff$ is
monotone non-decreasing in $[0,1]^N$, i.e.
\begin{equation}\label{lambda}
\exists\,\lambda>0:\quad\lambda|\bv-\bw|^2-\big(\bff(x,t,\bv)-\bff(x,t,\bw)\big)\cdot(\bv-\bw)\ge
0,\quad(x,t)\in\,Q,\ \forall\,\bv,\bw\in [0,1]^N,
\end{equation}
in particular, if $\bff$ is Lipschitz continuous in $\bv$, then
there exists at most one solution $\bu$ of the variational
inequality  (\ref{vit}) in the class  (\ref{regu}) and initial
condition $\bu_0\in\,\tilde{\K}$.

In order to prove the uniqueness of solution, we suppose that
$\bu_1$ and $\bu_2$ are two solutions of the variational inequality
(\ref{vit}) with initial condition $\bu_0\in\,\tilde{\K}$ and
$\bff=\bff(\bu_1)$, $\bff=\bff(\bu_2)$ respectively. Then, choosing
$\bu_2$ and $\bu_1$ as test functions, respectively, using
(\ref{lambda}) we find
\begin{multline*}
\frac12\int_\Omega|\bu_2(t)-\bu_1(t)|^2+\int_{Q_t}|\nabla(\bu_2-\bu_1)|^2\\
\le\int_{Q_t}\big(\bff(\bu_2)-\bff(\bu_1)\big)\cdot\big(\bu_2-\bu_1)\le\lambda\int_{Q_t}|\bu_2-\bu_1|^2
\end{multline*}
and so, by Grownvall inequality $\bu_1=\bu_2$ a.e. in $Q$, since
$\bu_1(0)=\bu_2(0)=\bu_0$.

We redefine the variational formulation of the approximating problem
(\ref{iveps})  in the framework of this section with $\Theta_\eps$
defined in (\ref{thetaeps2}) and  with initial condition only in
$L^2(\Omega)^N$,
\begin{equation}\label{ivepst}
\int_Q\partial_t\bu^\eps\cdot\bv+\int_Q\nabla\bu^\eps\cdot\nabla\bv+\int_Q\Theta_\eps(\bu^\eps)\cdot\bv
=\int_Q\bff(\bu^\eps)\cdot \bv,\quad\forall\bv\in\,
L^2(0,T;H^1(\Omega)^N)\cap L^\infty(Q)^N.
\end{equation}

Arguing as in Theorem \ref{fu} we may prove the existence of a
solution of the approximating problem (\ref{aprox}), with initial
condition $\bu_0\in\tilde{\K}$ as long as $\bff$
 satisfies (\ref{fnova1}) and  (\ref{fnova2}). We also have
uniqueness   if we
 assume (\ref{lambda}).

 \begin{thm}  Suppose that $\bff$ satisfies {\em (\ref{fnova1})},
{\em (\ref{fnova2})} and
 {\em (\ref{lambda})} and $\bu_0\in\tilde{\K}$.

  Let $\bu^\eps$ and $\bu$ be, respectively, the unique solution of the approximating
 problem {\em(\ref{aprox})} and of the variational inequality
 {\em(\ref{vit})}, both with initial condition $\bu_0$. Then there exists a
 positive constant $c=c(\varphi_1,T)$ such that  the
 following estimate in the norm {\em (\ref{norm})} of $\bW=C([0,T];L^2(\Omega)^N)\,\cap\,
 L^2(0,T;H^1(\Omega)^N)$ holds,
 \begin{equation}\label{eps}
\|\hspace{-1pt}|\bu^\eps-\bu\|\hspace{-1pt}|\le c\sqrt{\eps}.
  \end{equation}
 \end{thm}
 \proof We choose in (\ref{ivepst}) $\bv=\bu^\eps-\bu$ as test
 function. Since $\bu\in\bk$, then
 $$\int_Q\Theta_\eps(\bu^\eps)\cdot(\bu^\eps-\bu)\ge 0$$
 and so
 \begin{equation}\label{ueue}
\int_{Q_t}\partial_t\bu^\eps\cdot(\bu^\eps-\bu)+\int_{Q_t}\nabla\bu^\eps\cdot(\bu^\eps-\bu)\le\int_{Q_t}\bff(\bu^\eps)\cdot(\bu^\eps-\bu).
 \end{equation}

 Choosing, as test function in (\ref{vit})
 $\bv^\eps=\big((\uep_1-\frac\eps{N})^+,\ldots,(\uep_N-\frac\eps{N})^+)$
 we get
  \begin{multline}\label{uu}
\int_{Q_t}\partial_t\bu\cdot(\bu^\eps-\bu)+\int_{Q_t}\nabla\bu\cdot\nabla(\bu^\eps-\bu)\\\ge\int_{Q_t}\bff(\bu)\cdot(\bu^\eps-\bu)
+\int_{Q_t} \big[P\bu-\bff(\bu)\big]\cdot(\bu^\eps-\bv^\eps)
 \end{multline}
 and subtracting (\ref{uu})  from   (\ref{ueue})  we get
\begin{alignat}{1}\nonumber\frac12\int_{\Omega}|\bu^\eps(t)-\bu(t))|^2&+\int_{Q_t}|\nabla(\bu^\eps-\bu)|^2\\
\label{esteps}&\le\int_{Q_t}(\bff(\bu^\eps)-\bff(\bu))\cdot(\bu^\eps-\bu)+\int_{Q_t}
\big[P\bu-\bff(\bu)\big]\cdot(\bv^\eps-\bu^\eps)\\
\nonumber&\le\lambda\int_{Q_t}|\bu^\eps-\bu|^2+\eps\int_{Q_t}
|P\bu-\bff(\bu)|,
\end{alignat}
since $\|\bv^\eps-\bu^\eps\|_{L^\infty(Q)^N}\le \eps$. Letting
$C=C(\varphi_1,T)=\|P\bu-\bff(\bu)\|_{L^1(Q)^N}$ and dropping the
nonnegative term $\displaystyle{\int_{Q_t}|\nabla(\bu^\eps-\bu)|^2}$
in (\ref{esteps}) we obtain, by application of the Grownwall
inequality,
$$\int_\Omega|\bu^\eps(t)-\bu(t)|^2\le 2\eps C e^{2\lambda t}$$
and using again (\ref{esteps}), also
$$|\hspace{-1.5pt}\|\bu^\eps-\bu\|\hspace{-1.5pt}|\le c\sqrt{\eps}.$$
\fim\vspace{2mm}

With similar arguments we may give a continuous dependence result
for solutions of the variational inequality (\ref{vi2}).

Suppose we have a sequence $\tende{\bf\bff^\nu}{\nu}{\bff}$ in the
following sense
\begin{equation}\label{4646}
\left.\begin{array}{c} \bff^\nu=\bff^\nu(x,t,\bv)\mbox{ are
continuous in
}\bv\in\,[0,1]^N,\mbox{ for a.e. }(x,t)\in Q\vspace{3mm}\\
\tende{\bff^\nu(\cdot,\cdot,\bv)}{\nu}{\bff(\cdot,\cdot,\bv)}\mbox{
in }L^1(Q)^N\mbox{ for all fixed }\bv\in\,[0,1]^N.
\end{array}\right\}
\end{equation}

In addition, the assumption (\ref{fnova2}) is satisfied for all
$\bff$ uniformly in $\nu$, i.e., there is a common $\varphi_1$ such
that (\ref{fnova2}) holds for all $\nu$, and the initial data are
such that
\begin{equation}\label{4747}
\tilde{\K}\ni\tende{\bu_0^\nu}{\nu}{\bu_0}\qquad\mbox{ in
}L^2(\Omega)^N.
\end{equation}

Hence, by Theorem \ref{fu}, it is clear that there are solutions
$\{\bu^\nu\}_{\nu\in\N}$ to the corresponding problems associated
with $\bff^\nu$ and $\bu_0^\nu$ and, moreover, they satisfy
(\ref{regu}) and (\ref{regut}) uniformly in $\nu$, i.e., their norms
in those spaces are bounded by a constant independent of $\nu$.
Therefore, we have a function $\bu$ in the same class (\ref{regu})
and (\ref{regut}), and a subsequence, still denoted by $\nu$, such
that
\begin{equation}\label{erro}
\bu^\nu\,\underset\nu\lraup\,\bu\mbox{ in
}L^2(0,T;H^1(\Omega)^N)\mbox{ weak$\quad $and in$\quad$ }
L^\infty(0,T;\tilde{\K})\mbox{ weak-}* \end{equation}
\begin{equation}
\tende{\bu^\nu}{\nu}{\bu}\quad\mbox{ a.e. in }Q\quad\mbox{ and
in }\quad L^p(Q)^N,\quad\forall\,1\le p<\infty.
\end{equation}

By assumption (\ref{4646}) and Lebesgue Theorem, we conclude first
that $\displaystyle{\tende{\bff^\nu(\bu^\nu)}{\nu}{\bff(\bu)}}$ a.e.
in $Q$ and in $L^1(Q)^N$, as well as
\begin{equation}\label{fnu}
\tende{\displaystyle{\int_Q\bff^\nu(\bu^\nu)\cdot\bu^\nu}}{\nu}{\displaystyle{\int_Q\bff(\bu)\cdot\bu}},
\end{equation}
\begin{equation}\label{fnu0}
\tende{\displaystyle{\int_Q\big(\bff^\nu(\bu^\nu)-\bff(\bu)\big)
\cdot(\bu^\nu-\bu)}} {\nu}{0},
\end{equation}
since, in particular, $|\bu^\nu|\le 1$ and $|\bu|\le 1$ a.e. in $Q$.

Recalling (\ref{puiu}) for each $\nu$, we may take the limit in
\begin{equation}\label{pfnu}
\int_Q\big(P\bu^\nu-\bff^\nu(\bu^\nu)\big)\cdot(\bv-\bu^\nu)\ge 0
\end{equation}
for a fixed $\bv\in L^2(0,T;\tilde{\K})$. Using (\ref{fnu}) and
(\ref{erro}), that in particular imply
$$P\bu\in L^1(Q)^N\qquad\mbox{ and }\qquad\liminf_\nu\int_{Q_t}
P\bu^\nu\cdot\bu^\nu\ge\int_{Q_t}P\bu\cdot\bu,\quad\forall\,t\in(0,T),$$
we conclude that $\bu$ is a solution of (\ref{vi2}) with initial
condition $\bu_0$.

Using $\bv=\bu\cchi_{(0,t)}+\bu^\nu\cchi_{(t,T)}$ in (\ref{pfnu})
and $\bv=\bu^\nu\cchi_{(0,t)}+\bu\cchi_{(t,T)}$ in (\ref{vi2}) we
find, for a.e. $t\in(0,T)$,
$$\frac12\int_\Omega|\bu^\nu(t)-\bu(t)|^2+\int_{Q_t}|\nabla(\bu^\nu-\bu)|^2\le
\int_{Q_t}\big[\bff^\nu(\bu^\nu)-\bff(\bu)\big]\cdot(\bu^\nu-\bu)+
\frac12\int_\Omega|\bu_0^\nu-\bu_0|^2$$ and, by (\ref{fnu0}), we
conclude that $\displaystyle{\tende{\bu^\nu}{\nu}{\bu}}$ strongly in
$\bW$. Therefore, we have proved the following result

\begin{thm} \label{t34} If $\bu^\nu$ denotes the solution to the variational
inequality {\em (\ref{vi2})} with $\bff^\nu$ satisfying the
assumptions {\em (\ref{4646})} and {\em (\ref{fnova2})} uniformly in
$\nu$ and initial condition satisfying {\em (\ref{4747})}, then
there exists a subsequence $\{\bu^\nu\}_{\nu\in\N}$ such that
$$\tende{\bu^\nu}{\nu}{\bu}\quad\mbox{ in }\quad
C([0,T];L^2(\Omega)^N\cap\tilde{\K})\cap L^2(0,T;H^1(\Omega)^N)\cap
L^p(Q)^N,\quad\forall\,1\le  p<\infty,$$ where $\bu$ is a solution
to {\em (\ref{vi2})} corresponding to the limit $\bff$ and the limit
initial condition $\bu_0$. In addition, if $\bff$ satisfies {\em
(\ref{lambda})}, by uniqueness of $\bu$, the whole sequence
$\{\bu^\nu\}_{\nu\in\N}$ converges.
\end{thm}
\fim

\vspace{5mm}

\section{The multiphases system and its characterization} $ $

In this section we consider a variational solution $\bu$ of
(\ref{vit}) obtained in Theorem \ref{fu}, i.e., satisfying
(\ref{regu}) and  (\ref{regut}). Setting
\begin{equation} w_i(\bu)=1-\sum_{j\neq i}u_i,\qquad
i=1,\ldots,N,
\end{equation}
each component $u_i$ satisfies a double obstacle problem
\begin{equation}
0\le u_i(x,t)\le w_i(x,t)\qquad\mbox{ a.e. }(x,t)\in Q,\qquad
i=1,\ldots,N.
\end{equation}

For an arbitrary nonnegative and bounded function
$\varphi=\varphi(x,t)$ defined for $(x,t)\in Q$, such that
\begin{equation}
\K_0^\varphi=\{v\in L^2(0,T;H^1(\Omega):\ 0\le
v\le\varphi\ \mbox{ in }Q\}\neq\emptyset,
\end{equation}
and for a given $g\in L^1(Q)$, we may introduce the parabolic double
obstacle scalar problem
\begin{equation}\label{dop}
u\in\K_0^\varphi:\qquad\int_Q\partial_tu(v-u)+\int_Q\nabla
u\cdot\nabla(v-u)\ge \int_Q g(v-u)\qquad\forall v\in\K_0^\varphi,
\end{equation}
subject to a given compatible initial condition
\begin{equation}\label{45}
u(0)=u_0\qquad\mbox{ in }\Omega.
\end{equation}

For each $i=1,\ldots,N$, we have $u_i\in\K_0^{w_i}$ and, by choosing
in (\ref{vit}) $\bv\in L^2(0,T;\K)$, such that $v_j=u_j$ for $j\neq
i$ and $v_i=v\in\K_0^{w_i}$ arbitrarily, it is clear that $u_i$ is a
solution of the scalar double obstacle problem (\ref{dop}) with
$\varphi=w_i$ and $g=f_i(\bu)$. Hence we can obtain further
properties of our solution by applying the general theory of the
obstacle problem. For the sake of completeness we prove here the
result below.

Let
\begin{equation}\label{46}
\varphi\in L^2(0,T;H^1(\Omega))\cap L^\infty(Q)\ \mbox{ with
}\varphi\ge 0\mbox{ a.e. in }Q,
\end{equation}
\begin{equation}\label{47}
\partial_t\varphi\in
L^2\big(0,T;\big(H^1(\Omega)\big)'\big)\ \mbox{ with }P\varphi\in
L^1(Q),\qquad\frac{\partial\varphi}{\partial\bn}=0\mbox{ on
}\partial\Omega\times(0,T),
\end{equation}
and
\begin{equation}\label{48}
g\in L^1(Q),\qquad u_0\in L^2(\Omega),\quad 0\le
u_0\le\varphi(0)\ \mbox{ in }\Omega.
\end{equation}

We observe that (\ref{47}) means that $\varphi$ satisfies the
formula
$$\langle\partial_t\varphi,v\rangle_t=\int_{Q_t}v\,P\varphi
-\int_{Q_t}\nabla\varphi\cdot\nabla v,\qquad\forall
v\in\,L^2(0,T;H^1(\Omega))\cap L^\infty(Q).$$

\begin{prop} \label{p41} Under the assumptions {\em (\ref{46})-(\ref{48})} the
unique solution $u\in\K_0^\varphi$ to the scalar problem {\em
(\ref{dop})-(\ref{45})} is such that
\begin{equation}\label{49}
u\in\,C([0,T];L^2(\Omega))\,\cap\,L^\infty(Q),\qquad\partial_tu\in\,L^1(Q)+L^2(0,T;\big(H^1(\Omega)'\big),
\end{equation}
and it satisfies the parabolic semilinear equation
\begin{equation}\label{410}
Pu=g+g^-\cchi_{\{u=0\}}-(P\varphi-g)^-\cchi_{\{u=\varphi\}}\qquad\mbox{
a.e. in }Q.
\end{equation}
\end{prop}
\proof Using the function $\theta_\eps$ given by (\ref{thetaeps})
and defining
\begin{equation}\label{411}
\vartheta_\eps(v)=g^-\theta_\eps(v)-(P\varphi-g)^-\theta_\eps(\varphi-v)
\end{equation}
we can consider the approximating problem, for $\eps>0$,
\begin{equation}\label{412}
\int_Q\big(P\uep+\vartheta_\eps(\uep)\big)v=\int_Qgv,\qquad\forall\,v\in\,L^2(0,T;H^1(\Omega))\,\cap\,L^\infty(Q),
\end{equation}
with the initial condition $\uep(0)=u_0$ in $\Omega$. Since
$\vartheta_\eps$ is monotone and $\varphi$ is bounded, arguing as in
Theorem \ref{fu}, the problem (\ref{412}) has a unique solution
$\uep$ in the class (\ref{49}). Moreover, it satisfies
\begin{equation}
-\eps\le\uep\le\varphi+\eps\qquad\mbox{ a.e. in }Q,
\end{equation}
as we can show by choosing, in (\ref{412}), $v=(-\uep-\eps)^+$ and
$v=(\uep-\varphi-\eps)^+$, respectively. Indeed, in the first case
we have
$$\int_QvPv=-\int_QvP\uep=\int_{\{v>0\}}v\left(\vartheta(\uep)-g\right)=\int_{\{\uep<-\eps\}}(-g^--g)\le
0,$$ since $\vartheta_\eps(\uep)=-1$ and
$\vartheta_\eps(\varphi-u_\eps)=0$, because $\uep<-\eps$ and
$\varphi-u_\eps>\eps$, and, in the second case,
\begin{alignat*}{1}\int_QvPv=\int_QvP(\uep-\varphi)=&\int_{\{v>0\}}v\left(g-\vartheta(\uep)-P\varphi\right)\\
=&\int_{\{\varphi-\uep>\eps\}}
\left(-(P\varphi-g)-(P\varphi-g)^-\right)\le 0,\end{alignat*} since
$\vartheta_\eps( \varphi-\uep)=-1$ and $\vartheta_\eps(\uep)=0$ if
$\varphi-u_\eps<-\eps$ and $\uep>\varphi+\eps$.

Hence, using the monotonicity argument, we easily conclude that
$u=\displaystyle{\lim_{\eps\rightarrow 0}\uep}\in\K_0^\varphi$ is
the unique solution of the variational inequality (\ref{dop}).
Remarking that, from (\ref{411}) we have
$$-g^-\le\vartheta_\eps(\uep)\le(P\varphi-g)^-\qquad\mbox{ a.e. in
}Q,$$  from (\ref{412}) we deduce in the limit the Lewy-Stampacchia
inequalities
$$(P\varphi-g)^-\le Pu-g\le g^-\qquad\mbox{ a.e. in }Q.$$

In particular, this yields $Pu\in\,L^1(Q)$ and (\ref{dop}) implies
that $u$ also solves
\begin{equation}\label{414}
\int_Q\big(Pu-g\big)(v-u)\ge
0,\qquad\forall\,v\in\tilde{\K}_0^\varphi,
\end{equation}
where $\tilde{\K}_0^\varphi=\{v\in\,L^2(Q):\,0\le v\le\varphi\mbox{
in }Q\}\,\subset\,L^\infty(Q)$.

Let $\mathcal{O}\subset Q$ be an arbitrary measurable set and set
$v=u$ in $Q\setminus\mathcal{O}$ and $v=\delta\varphi$ in
$\mathcal{O}$, with $\delta\in\,[0,1]$, in (\ref{414}). Since
$\mathcal{O}$ is arbitrary, we conclude the pointwise inequality
\begin{equation}\label{416}
\big(Pu-g\big)(\phi-u)\ge 0\qquad\forall\,\phi\in\,[0,\varphi(x,t)]\
\mbox{ a.e. in }Q,
\end{equation}
which implies, up to null measure subsets of $Q$,
\begin{equation}\label{417}
Pu-g\ge 0\ \mbox{ in }\ \{u=0\},\qquad Pu-g\le 0\ \mbox{ in }\
\{u=\varphi\},
\end{equation}
\begin{equation}
Pu=g\ \mbox{ in }\ \Lambda=\{0<u<\varphi\}.
\end{equation}

On the other hand, arguing as in Lemma 2 of \cite{r3} and noting
that $V=(u,-\nabla u)\in\,L^1(Q)^{n+1}$ and $D\cdot V=Pu\in\,
L^1(Q)$, with $D=(\partial_t,\partial_{x_1},\ldots,\partial_{x_n})$,
we have
$$Pu=0\ \mbox{ a.e. in }\ \{u=0\}\qquad\mbox{ and }\qquad
Pu=P\varphi\ \mbox{ a.e. in }\ \{u=\varphi\}.$$

Hence, by (\ref{416}), up to neglectable sets, we have
$\{u=0\}\subset\{g\le 0\}$ and $\{u=\varphi\}\subset\{P\varphi\le
g\}$, and using also (\ref{417}), we finally conclude
(\ref{410}).\fim\vspace{2mm}

\begin{thm} Any solutions $\bu$ of the variational inequality
{\em (\ref{vit})} {\em (}or {\em (\ref{vi2})}{\em )} under the
conditions of Theorem {\em \ref{fu}} satisfy the semilinear
parabolic system
\begin{multline}\label{418}
Pu_i=f_i(\bu)+f_i^-(\bu)\cchi_{\{u_i=0\}} \\
-\soma \frac1k(f_{i_1}(\bu)+\cdots+f_{i_k}(\bu))^+\cchi_{i_1\ldots
i_k}\qquad\mbox{ a.e. in }Q,
\end{multline}
where $\cchi_{i_1\ldots i_k}=\cchi_{I_{i_1\ldots i_k}},$ for
$k=1,\ldots,N$, denotes the characteristic function of
\begin{equation}
\label{ii1} I_{i_1\ldots i_k}=\{(x,t)\in Q:\,u_{i_1\ldots
i_k}(x,t)=1,\, u_{i_j}(x,t)>0\mbox{ for all } j=1,\ldots,k\}.
\end{equation}
\end{thm}
\proof We notice that the regularity (\ref{46}), (\ref{47}), holds
for $\displaystyle{w_i=1-\sum_{j\neq i}u_j}$, so $w_i$ can be chosen
as the upper obstacle of each component $u_i$, $i=1,\ldots,N$, of
$\bu$, to which we can apply the conclusions of Proposition
\ref{p41}. Since $\{u_i=0\}\subset\{f_i(\bu)\le 0\}$ a.e., for each
$i=1,\ldots,N$, we have
\begin{equation}\label{coinc}
Pu_i=f_i(\bu)+f_i^-(\bu)\cchi_{\{u_i=0\}}-(Pw_i-f_i(\bu))^-\cchi_{\{u_i=w_i,\,
u_i> 0\}} \qquad\mbox{ in }Q,
\end{equation}
and the condition (\ref{418}) will follow if we show that
\begin{equation}\label{iq}
(Pw_i-f_i(\bu))^-\cchi_{\{u_i=w_i,\, u_i> 0\}}=\soma
\frac1k\big(f_{i_1}(\bu)+\cdots+f_{i_k}(\bu)\big)^+\cchi_{i_1\ldots
i_k}\qquad\mbox{ in }Q,
\end{equation}

Observe that
\begin{equation*}
\big\{u_i=w_i,\, u_i>0\big\}=\uniao I_{i_1\ldots i_k},
\end{equation*} and these sets  are a.e. disjoint. Here the union is
taken also over all the subsets $\{i_1,\ldots,i_k\}$ of
$\{1,\ldots,N\}$ that include $i$ and over all $k=1,\ldots,N$. We remark that $Pw_i=Pu_i$ in
that subset and
\begin{itemize}
\item in the sets $I_i=\{u_i=1\}$, $Pw_i=0$ and
$(Pw_i-f_i(\bu))^-=f_i(\bu)^+$, for $i=1,\ldots,N$;
\item in each set $I_{i_1\ldots i_k}$, for $k\ge 2$, as we shall see,
$$(Pu_i- f_i(\bu))^-=\frac1k\big(f_{i_1}(\bu)+\cdots+f_{i_k}(\bu)\big)^+,$$  and this fact
concludes the proof.
\end{itemize}

Let $(x_0,t_0)\in I_{i_1\ldots i_k}$.  Recall that
$\{i_1,\ldots,i_k\}$ is the set of indexes for which we have
$0<u_{i_j}(x_0,t_0)$ (notice that $i\in\,\{i_1,\ldots,i_k\}$).
Denoting $\alpha=\min\{u_{i_j}(x_0,t_0):\,j=1,\ldots,k\}$, the set
$\mathcal{O}=\bigcap_{j=1}^k\{u_{i_j}>\alpha/2\}$ is measurable and
contains $(x_0,t_0)$. Given any measurable set $\omega\subset
\mathcal{O}$,
 choose, in (\ref{vi2}), as test function
$\bv=(v_1,\ldots,v_N)$ defined by
$$
v_{i_1}=u_{i_1}\pm\delta\cchi_\omega,\quad
v_{i_j}=u_{i_j}\mp\delta\cchi_\omega\mbox{ for a fixed
}j\in\{2,\ldots,k\},\quad  v_l=u_l\ \forall\, l\neq i_1,i_j,$$
 observing that
$$\sum_{j=1}^Nv_j=\sum_{j=1}^Nu_j\pm\delta\cchi_\omega\mp\delta\cchi_\omega=\sum_{j=1}^Nu_j\le
1$$ and
$$v_j\ge 0,\ j=1,\ldots,N,\qquad\mbox{ as long as }0<\delta\le\alpha/2.$$

Returning to the inequality (\ref{vi2}) and setting
$S_j=Pu_j-f_j(\bu)$, we get
$$\pm\delta\int_Q S_{i_1}\cchi_\omega\mp\delta\int_Q S_{i_j}\cchi_\omega\ge
0.$$

Since $\omega\supset\{(x_0,t_0)\}$ was taken arbitrarily in
$\mathcal{O}$ and $(x_0,t_0)$ is a generic point of $I_{i_1\ldots
i_k}$, we conclude that
\begin{equation}\label{qij}
S_{i_1}=S_{i_j},\qquad\mbox{ a.e.  in }I_{i_1\ldots
i_k,}\qquad\mbox{ for any }j\in\{2,\ldots,k\}.
\end{equation}

Recalling that $\displaystyle{\sum_{j=1}^N} Pu_j=Pu_{i_1\ldots
i_k}=0$, in the set $I_{i_1\ldots i_k}$ we get, using (\ref{qij}),
that
$$kS_{i_1}=S_{i_1}+\cdots+S_{i_k}=\big(Pu_{i_1}-f_{i_1}\big)+\cdots+\big(Pu_{i_k}-f_{i_k}\big)=Pu_1+\cdots+Pu_N-(f_{i_1}+\cdots+f_{i_k}),$$
where, for simplicity, we set $f_j=f_j(\bu)$, and so
$$S_i=S_{i_1}=-\frac1k\big(f_{i_1}+\cdots+f_{i_k}\big).$$

But in $I_{i_1\ldots i_k}$ we have $S_i\le 0$ (recall that $u_i=w_i$
and (\ref{417})) and so
\begin{equation*}(Pu_i-
f_i(\bu))^-=-(Pu_i-f_i(\bu))
=-S_i=\frac1k\big(f_{i_1}+\cdots+f_{i_k}\big)=\frac1k\big(f_{i_1}+\cdots+f_{i_k}\big)^+.
\end{equation*}
\fim\vspace{2mm}

\begin{cor}
Let $\bu$ be the solution of the variational inequality {\em
(\ref{vit})} {\em (}or {\em (\ref{vi2})}{\em )} under the conditions
of Theorem {\em \ref{fu}}.

Then, denoting by $|A|$ the $(n+1)$-Lebesgue measure of $A\subset
Q$, we have
\begin{equation}\label{contacto}
\Big|\big\{\,\sum_{j=1}^k
f_{i_j}(\bu)<0\big\}\,\cap\,\big\{\,\sum_{j=1}^k u_{i_j}=1,\,
u_{i_j}>0,\, j=1,\ldots,k\big\}\Big|=0
\end{equation}
for each partial coincidence subset $I_{i_1\ldots i_k}$, as well as
\begin{equation}\label{ls0}
\big|\{f_i(\bu)>0\}\cap\{u_i=0\}\big|=0,\qquad i=1,\ldots,N.
\end{equation}
\end{cor}
\proof Being $I_{i_1\ldots i_k}$ defined in (\ref{ii1}), using the
equation (\ref{418}), we obtain, for each $i_j$ with $j=1,\ldots,k$,
denoting $f_{i_j}=f_{i_j}(\bu)$,
$$Pu_{i_j}=f_{i_j}-\frac1k\big(f_{i_1}+\cdots+f_{i_k}\big)^+\qquad\mbox{ a.e in }I_{i_1\ldots i_k}.$$

Summing these  $k$ equations, we have
$$0=\sum_{j=1}^kPu_{i_j}=f_{i_1}+\cdots+f_{i_k}-\big(f_{i_1}+\cdots+f_{i_k}\big)^+=\big(f_{i_1}+\cdots+f_{i_k}\big)^-
\qquad\mbox{ a.e in }I_{i_1\ldots i_k}.$$

So, in $I_{i_1\ldots i_k}=\big\{\,\sum_{j=1}^k u_{i_j}=1,\,
u_{i_j}>0,\, j=1,\ldots,k\big\}$ we have
$\displaystyle{\sum_{j=1}^kf_{i_j}\ge 0}$ a.e. and (\ref{contacto})
follows.

The proof of (\ref{ls0}) is similar (recall (\ref{417})).

\fim\vspace{2mm}

As a consequence of this corollary the semilinear system (\ref{418})
can, in fact, be written in the equivalent form for $i=1,\ldots,N$,
\begin{multline}\label{eq70}
Pu_i=f_i(\bu)-f_i(\bu)\cchi_{\{u_i=0\}} \\
-\soma \frac1k(f_{i_1}(\bu)+\cdots+f_{i_k}(\bu))\cchi_{i_1\ldots
i_k}\qquad\mbox{ a.e. in }Q,
\end{multline}
since $\{u_i=0\}\subset\{f_i(\bu)\le 0\}$ and
$\displaystyle{I_{i_1\ldots
i_k}\subset\big\{\sum_{j=1}^kf_{i_j}(\bu)\ge 0\big\}}$ up to a
neglectable subset of $Q$.

This remark combined with the continuous dependence of the
variational solutions obtained in Theorem \ref{t34} yields an
interesting criteria of local stability of the characteristic
functions of the coincidence sets in the Lebesgue measure. Denote
$$\cchi^\nu_{i_1\ldots i_k}=\cchi_{\{u_{i_1\ldots i_k}^\nu=1,\
u_{i_j}^\nu>0\ \forall j=1,\ldots,k\}},\qquad \mbox{\small $1\le
i_1<\cdots<i_k\le N,\ i\in\{i_1,\ldots,i_k\}$}.$$

\begin{thm} Let the assumptions and notations of {\em Theorem
\ref{t34}} hold. Suppose that in some subset of positive measure
$\omega\subseteq Q$ the following assumption on the limit  problem
holds \begin{equation} \label{festab} \sum_{j=1}^kf_{i_j}(\bu)\neq
0\qquad\mbox{ a.e. in }\omega,\qquad 1\le i_1<\cdots<i_k\le N,\quad
k=1,\ldots,N.
\end{equation}

Then the associated characteristic functions are such that
\begin{equation}
\tende{ \cchi_{\{u_i^\nu=0\}} }{\nu}{\cchi_{\{u_i=0\}}}\qquad\mbox{
in }L^p(\omega),\qquad \forall\,i=1,\ldots,N,
\end{equation}
\begin{equation}
\tende{\cchi_{i_1\ldots i_k}^\nu}{\nu}{\cchi_{i_1\ldots
i_k}}\qquad\mbox{ in }L^p(\omega),\qquad \forall\,i_1,\ldots ,i_k,
\end{equation}
for all $p$, $1<p<\infty$.
\end{thm}
\proof We observe that each $\bu^\nu$   solves the system
\begin{equation}\label{system2}
Pu_i^\nu=f_i^\nu-f_i^\nu\cchi_{\{u_i^\nu=0\}}- \soma
\frac1k(f_{i_1}^\nu+\cdots+f_{i_k}^\nu)\,\cchi_{i_1\ldots
i_k}^\nu\qquad\mbox{ a.e. in }Q
\end{equation}
where, for simplicity, we set $f_i^\nu=f_i^\nu(\bu^\nu)$. By the
convergence $\displaystyle{\tende{\bu^\nu}{\nu}{\bu}}$, we have
$\tendef{P\bu^\nu}{\nu}{P\bu}$ in the distributional sense. Since
$0\le\cchi_{i_1\ldots i_k}^\nu\le 1$, there exists $\cchi_{i_1\ldots
i_k}^*$, with $0\le \cchi_{i_1\ldots i_k}^*\le 1$ in $Q$, such that
$$\tendef{\cchi_{i_1\ldots i_k}^\nu}{\nu}{\cchi_{i_1\ldots
i_k}^*}\qquad\mbox{ in }\quad L^\infty(Q)\mbox{ weak-}*.$$

Analogously, for some $\cchi_{i,0}^*$, with $0\le \cchi_{i,0}^*\le
1$ in $Q$, we have
$$\tendef{\cchi_{\{u_i^\nu=0\}}}{\nu}{\cchi_{i,0}^*}\qquad\mbox{ in }\quad L^\infty(Q)\mbox{ weak-}*.$$

We are going to prove that, in fact,
$$\cchi_{i,0}^*
=\cchi_{\{u_i=0\}}\qquad\mbox{  and }\qquad \cchi_{i_1\ldots
i_k}^*=\cchi_{i_1\ldots i_k}\quad\mbox{ a.e. in }\omega,$$which
concludes the proof, since the weak convergence to characteristic
functions in $L^p(\omega)$ is in fact strong, as it is well known.

Passing to the limit in (\ref{system2}), we obtain
$$Pu_i=f_i-f_i\cchi^*_{i,0}-
\soma \frac1k(f_{i_1}+\cdots+f_{i_k})\,\cchi_{i_1\ldots
i_k}^*\qquad\mbox{ a.e. in }Q$$ where, for simplicity, we have also
set $f_{i_j}=f_{i_j}(\bu)$.

But each $u_i$ also solves the equation (\ref{eq70}), so, by
subtraction, we obtain a.e. in $Q$,

\begin{equation}\label{igualdade}-f_i\big(\cchi_{\{u_i=0\}}-
\cchi^*_{i,0}\big)
-\soma\frac1k(f_{i_1}+\cdots+f_{i_k})\big(\cchi_{i_1\ldots
i_k}-\cchi_{i_1\ldots i_k}^*\big)=0.\end{equation}

Noticing that $\cchi_{\{u_i^\nu=0\}}u_i^\nu=0,$ passing to the
limit, we get $\cchi_{i,0}^*u_i= 0,$ which means that
$\cchi_{i,0}^*=0$ whenever $u_i>0$. To conclude that
$\cchi_{i,0}^*=\cchi_{\{u_i=0\}}$ we only need to prove that
 $\cchi_{i,0}^*=1$ if $u_i=0$.

 Recall that the sets
${\{u_i=0\}}$ and $I_{i_1\ldots i_k}$, $1\le i_1<\ldots<i_k\le N$,
$i\in\{i_1,\ldots,i_k\}$, $k=1,\ldots,N$, are mutually  disjoint.
Hence in $\{u_i=0\}$ we obtain
$$-f_i(1-\cchi^*_{i,0})+\soma\frac1k(f_{i_1}+\cdots+f_{i_k})\cchi_{i_1\ldots
i_k}^*=0$$ and since the left hand side is nonnegative, by the
assumption (\ref{festab}) we conclude that
$$\cchi^*_{i,0}=1\qquad\mbox{ and }\qquad\cchi_{i_1\ldots i_k}^*=0\qquad\mbox{
in }\{u_i=0\}\cap\omega.$$

Since $\cchi_{i_1\ldots i_k}^\nu(1-u_{i_1\ldots i_k}^\nu)= 0$ a.e.
in $Q$, taking the limit in $\nu$, we also obtain\linebreak
$\cchi_{i_1\ldots i_k}^*(1-u_{i_1\ldots i_k})= 0$ a.e in $Q$, i.e.
$\cchi_{i_1\ldots i_k}^*=0$ if $u_{i_1\ldots i_k}<1$. It remains to
evaluate $\cchi_{i_1\ldots i_k}^*$ when $u_{i_1\ldots i_k}=1$ and
$u_{i_j}>0$, for all $j=1,\ldots,k$ or when $u_{i_j}=0$, for some
$j=1,\ldots,k$.

In this later case, where $u_{i_j}=0$, for some $j=1,\ldots,k$, we
have $\cchi_{i_1\ldots i_k}=0$ and, since we already know that
$\cchi_{\{u_{i_j}=0\}}=\cchi_{i_j,0}^*$, from (\ref{igualdade}) for
the index $i_j$, we get
$$
\displaystyle{\sum_{\mbox{\tiny $\begin{array}{c}
1\!\!\le i_1\!\!\!<\!\!\ldots\!\!<i_k\!\!\!\le N\vspace{1mm}\\
i_j\in\{i_1,\ldots,i_k\}\end{array}$}}}
\frac1k(f_{i_1}+\cdots+f_{i_k})\cchi_{i_1\ldots i_k}^*=0.$$

Then, by the assumption (\ref{festab}) we have $\cchi_{i_1\ldots
i_k}^*=0$ in $\big(Q\setminus I_{i_1\ldots i_k}\big)\cap\omega$.

Finally, in $I_{i_1\ldots i_k}\cap\omega$, again from
(\ref{igualdade}), we obtain
$$\frac1k(f_{i_1}+\cdots+f_{i_k})\big(1-\cchi_{i_1\ldots i_k}^*\big)=0$$
and the assumption (\ref{festab})  yields that $\cchi_{i_1\ldots
i_k}^*=1$, completing the proof. \fim

\vspace{5mm}

Acknowledgment: Research partially supported by FCT Project POCI/MAT/57546/2004.

\vspace{5mm}

\begin{tabular}{ll}
Jos\'{e} Francisco Rodrigues$\qquad\qquad$&$\qquad\qquad$Lisa Santos\\
Universidade de Lisboa/CMAF$\qquad\qquad$&$\qquad\qquad$Universidade do Minho/CMat\\
Prof. Gama Pinto 2$\qquad\qquad$&$\qquad\qquad$Campus de Gualtar\\
 1649--003 Lisboa$\qquad\qquad$&$\qquad\qquad$4710--057 Braga\\
  Portugal$\qquad\qquad$&$\qquad\qquad$Portugal\\
  $\qquad\qquad$&$\qquad\qquad$\\
rodrigue@fc.ul.pt$\qquad\qquad$&$\qquad\qquad$lisa@math.uminho.pt
\end{tabular}

\end{document}